\documentclass[12pt]{article}

\usepackage{amsmath,amsthm,graphicx}
\title{A note on the history of the four-colour conjecture}

\author{Brendan D. McKay\\
\small Research School of Computer Science\\[-0.8ex]
\small Australian National University\\[-0.8ex]
\small Canberra, ACT 0200, Australia\\
\small\tt bdm@cs.anu.edu.au}

\date{}

\begin{document}
\maketitle

\begin{abstract}
 \noindent
 The four-colour conjecture was brought to
 public attention in 1854, most probably by Francis or Frederick Guthrie.
 This moves back by six years the date of the earliest known publication.
\end{abstract}

As the most famous problem in discrete mathematics, the four-colour
problem has been the subject of much historical as well as mathematical
investigation. We refer to the surveys of this history in two excellent recent
books~\cite{Fritsch,Wilson}.

As the early history is understood, Francis Guthrie communicated the problem
to his brother Frederick, who in turn showed his teacher Augustus De Morgan.
In what remains the earliest known document, De Morgan
wrote of the problem to Hamilton on the same day that he heard it,
October 23, 1852.  Some additional private correspondence is known from
the following few years, but it has been believed that the first publication of the
problem was in an anonymous book review by De Morgan in
the magazine \textit{The Athen{\ae}um} in 1860~\cite{deMorgan}.

However, there was an earlier appearance in the same magazine, apparently
overlooked until now. On June 10, 1854, a letter appeared in the
Miscellanea section of \textit{The Athen{\ae}um}~\cite{FG}.  Printed in smaller
type than the surrounding text, it read in full: 

\begin{quotation}\noindent
   \textit{Tinting Maps}.---In tinting maps, it is desirable for the sake of
   distinctness to use as few colours as possible, and at the same
   time no two conterminous divisions ought to be tinted the same.
   Now, I have found by experience that \textit{four} colours are necessary
   and sufficient for this purpose,---but I cannot prove that this is
   the case, unless the whole number of divisions does not exceed five.
   I should like to see (or know where I can find) a general proof of
   this apparently simple proposition, which I am surprised never
   to have met with in any mathematical work.\hfill F.\,G.~~
\end{quotation}

The magazine does not identify ``F.\,G.'', but the short period of time between
this letter and the known interaction between Francis and Frederick Guthrie
makes it highly likely that one of them was responsible.  It doesn't seem
possible to identify which of the brothers it was, but I favour Francis for
the following, inconclusive, reason. In 1880, Frederick carefully
attributed the discovery to Francis and did not mention having studied
the problem himself~\cite{PRSE}.  This does not seem consistent with
him having written a 1854 letter that takes all the credit for~it.

An intriguing, but highly conjectural, possibility is that ``F.\,G.'' was Francis
Galton.  At the time, Galton was famous for his explorations in Africa, and
later he wrote about the art of map making~\cite{Crilly}.  However, between the
1854 letter of ``F.\,G.'' and the earliest proven connection between Galton and
the four-colour problem lie more than 20 years~\cite{Crilly}.

The letter is also of interest in that ``F.\,G.'' admits he can't solve the
problem.  This is at variance with Frederick's 1880 testimony that Francis
had a proof in 1852 even though it ``did not seem altogether satisfactory
to himself''~\cite{PRSE}.  As we know, De Morgan in his 1860 article
believed
he was presenting sufficient justification to his readers, while in reality
he was only observing that $K_5$ is non-planar.

The library of City University London holds a ``marked copy'' of
\textit{The Athen{\ae}um} on which the editors wrote the identities of many
of the anonymous contributors~\cite{Marked}. Unfortunately, the identity
of ``F.\,G.'' is not indicated there.  However the author or editor
of the Miscellanea
section that includes the letter of ``F.\,G.'' is indicated. The unclear
handwriting probably reads ``Winton'', who remains unidentified.

Robin Wilson made useful comments on the first draft of this note.
Thanks also to Micheline Beaulieu, Norman Biggs, Colin Burrow,
Tony Crilly, Keith Lloyd, and Sheila Munton.

\begin{figure}
\centering
 \includegraphics[scale=0.8]{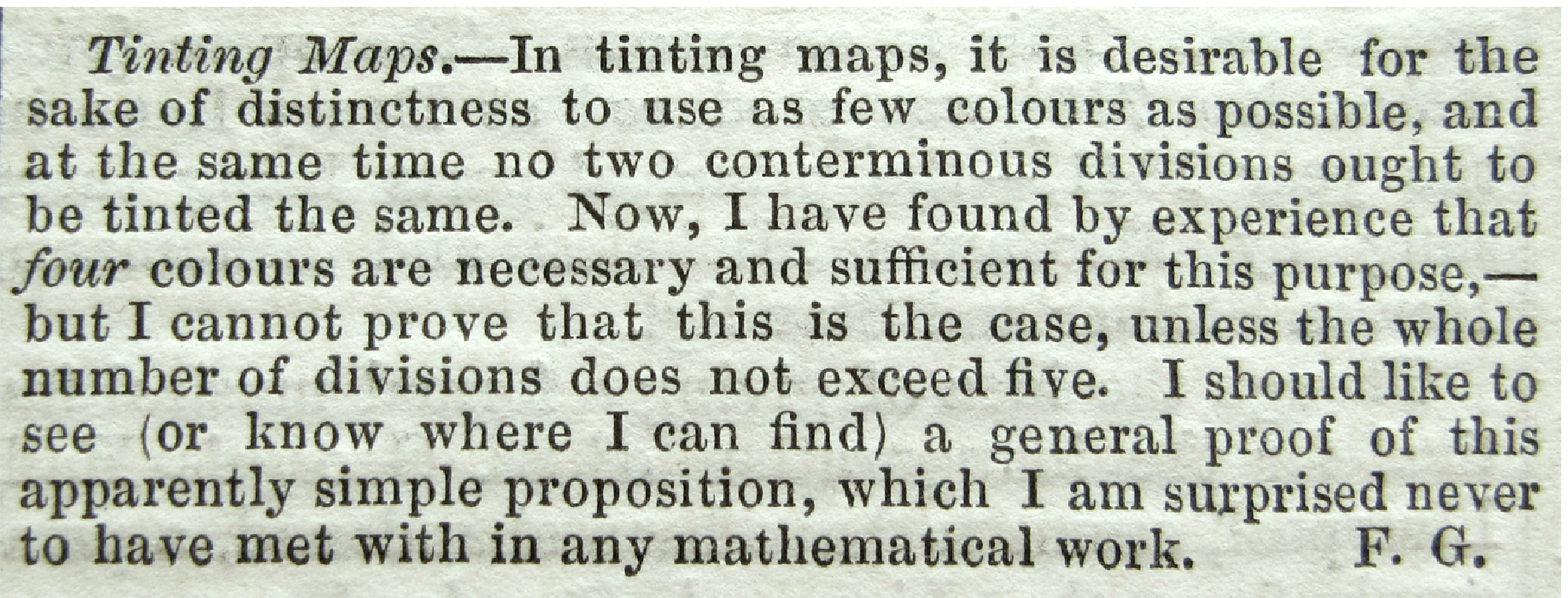}
 \caption{Letter in \textit{The Athen\ae um} of June 10, 1854.}
\end{figure}

\end{document}